\documentclass[11pt]{article}

\usepackage{amsmath}
\usepackage{amsfonts}

\newcommand{\be}{\begin}
\newcommand{\Ref}[1]{(\ref{#1})}

\newcommand{\cm}{\mathcal M}
\newcommand{\ca}{\mathcal A}
\newcommand{\cb}{\mathcal B}
\newcommand{\cg}{\mathcal G}

\def\fl{\mathcal R}

\newcommand{\al}{\alpha}
\newcommand{\del}{\delta}
\newcommand{\eps}{\epsilon}
\def\ga{\gamma}
\def\Ga{\Gamma}

\def\lla{\lambda}
\def\La{\Lambda}
\def\om{\omega}

\newcommand{\si}{\sigma}
\def\Si{\Sigma}

\newcommand{\co}{\mathbb C}

\def\e{{\mathbb E}}
\def\f{{\mathbb F}}
\def\bg{\mathbb G}
\def\bh{\mathbb H}

\def\bp{\mathbb P}
\def\q{{\mathbb Q}}
\newcommand{\R}{\mathbb R}

\newcommand{\z}{\mathbb Z}

\def\cp#1{{\mathbb{CP}^#1}}
\newcommand{\cpb}{\overline{\mathbb C\mathbb P}\vphantom{P}^2}

\newcommand{\uset}{\underset}
\newcommand{\oset}{\overset}
\newcommand{\uline}{\underline}
\newcommand{\oline}{\overline}

\newcommand{\la}{\langle}
\newcommand{\ra}{\rangle}
\newcommand{\st}{\;|\;}
\newcommand{\ti}{\tilde}
\newcommand{\prtl}{\partial}
\def\square{\kern20pt{\vbox{\hrule height.4pt
        \hbox{\vrule width.4pt height 6pt\kern6pt
                \vrule width.4pt}
        \hrule height.4pt}}}

\newcommand{\rpy}{\R_+\times Y}

\def\coker{\text{coker}}
\newcommand{\fact}{\frac1{8\pi^2}}
\def\inv{^{-1}}
\def\hol{\text{Hol}}

\def\proof{{\em Proof.}\ }
\def\loc{\text{loc}}

\def\frakg{\mathfrak g}

\def\U#1{\text{U}(#1)}
\def\u#1{\text{u}(#1)}
\def\SL{\text{SL}}
\def\SU#1{\text{SU}(#1)}
\def\su#1{\text{su}(#1)}
\def\SO#1{\text{SO}(#1)}

\def\Sp#1{\text{Sp}(#1)}
\def\Spin#1{\text{Spin}(#1)}
\def\hf{\text{HF}}

\def\thf{\widetilde{HF}\vphantom{F}}

\def\tsum{\textstyle\sum}
\def\tprod{\textstyle\prod}
\def\tcm{\tilde\cm}
\def\hcm{\widehat\cm}
\def\harm{\mathcal H}
\def\cl{\mathcal L}
\def\ct{\mathcal T}
\def\bt{t}
\def\bff{\mathbf F}
\def\tt{_{T,\bt}}
\def\torsion{\text{torsion}}

\def\cf{\text{CF}}
\def\lls#1#2{{L^{#1,\si}_#2}}
\def\ls#1{{L^{#1,\si}}}
\def\rbp{\R^{b^+}}

\def\iu{_{i\nu}}
\def\dsp{\vspace{.2cm}}
\def\dett{\text{det}}
\def\llw#1#2#3{{L^{#1,#2}_#3}}

\def\exot{e}

\begin{document}

\title{An inequality for the $h$-invariant\\
in instanton Floer theory}

\author{Kim A.\ Fr\o yshov\thanks{Partially supported by
a post-doctoral grant from the Norwegian Research Council,
and by NSF grant DMS-9971731.}}

\date{4 March 2003}

\maketitle

\bibliographystyle{plain}

\be{abstract}

Given a smooth, compact, oriented $4$-manifold $X$ with 
a homology sphere $Y$ as boundary and $b_2^+(X)=1$,
and given an embedded surface $\Si\subset X$
of self-intersection~$1$, we prove an inequality relating $h(Y)$,
the genus of $\Si$, and a certain invariant of
the orthogonal complement of $[\Si]$
in the intersection form of $X$.
\end{abstract}

\section{Introduction}

In \cite{Fr3} we defined a surjective group homomorphism from the
homology cobordism group of oriented (integral) homology $3$-spheres
to the integers,
\[h:\theta^H_3\to\z,\]
using instanton Floer theory.
The main purpose of this paper is to establish
the following property of this invariant.

\be{thm}\label{hthm}
Let $X$ be a smooth, compact, oriented $4$-manifold with a homology sphere
$Y$ as boundary and with $b_2^+(X)=1$.
Let $\Si\subset X$ be a closed surface
of genus $g$ and self-intersection $\Si\cdot\Si=1$, and 
let $\cl\subset H^2(X;\z)/\torsion$ be the
sublattice consisting of all vectors which
vanish on $[\Si]$. Then the intersection form of $X$ restricts to a unimodular
negative definite form on $\cl$, and we have
\[h(Y)+\lceil g/2\rceil\ge\exot(\cl),\]
where $\exot(\cl)$ is a non-negative integer which depends only on $\cl$.
Moreover,
\be{description}
\item[$\bullet$]$\exot(\cl)=0$ precisely when $\cl$ is diagonal
\item[$\bullet$]$\exot(-kE_8\oplus \cl)=k$ if $\cl$ is diagonal and $k\ge0$.
\end{description}
\end{thm}
By ``surface in $X$'' we mean a $2$-dimensional compact, oriented, connected
smooth submanifold of $X$.
$\lceil x\rceil$ denotes the smallest integer $\ge x$.

See Section~\ref{main-thm} for a slightly more general statement and the
definition of $\exot(\cl)$.
To compute the term $\lceil g/2\rceil$
in the inequality we rely on Mu\~noz'\ description \cite{Mun1}
of the ring structure of the Floer cohomology of the $\SO3$ bundle
$E\to S^1\times\Si$ where $E$ is the pull-back of the non-trivial 
$\SO3$ bundle over $\Si$.

Note that if $Z$ is a negative definite $4$-manifold with a homology
sphere $Y$ as boundary then one can apply the theorem 
to $X=Z\#\cp2$ with $g=0$.
When $Y=S^3$ one recovers Donaldson's diagonalization
theorem~\cite{D1,D2}.

As another example, we obtain bounds on how much $h$ may change under $\pm1$
surgery on knots:
\be{cor}\label{genus-bound}
Let $Y$ be an oriented homology $3$-sphere and $\ga$ a knot in $Y$
of slice genus $g$. If $Y_{\ga,-1}$ is the result of $-1$ surgery on $\ga$ 
then 
\[0\le h(Y_{\ga,-1})-h(Y)\le\lceil g/2\rceil.\]
\end{cor}

Here the slice genus may be defined as the smallest non-negative integer
$g$ for which there exists a smooth rational homology cobordism $W$
from $Y$ to some rational homology sphere $Y'$ and a genus $g$
surface $\Si\subset W$ such that $\prtl\Si=\ga$.
(We do not know whether this definition agrees with the usual one for $Y=S^3$.)

To deduce the corollary from the theorem, let $Z$ be the surgery cobordism
from  $Y$ to $Y_{\ga,-1}$ and set $W'=W\cup_{Y'}\oline W\cup_YZ$, which
is a cobordism from $Y$ to $Y_{\ga,-1}$. By attaching a 
suitable $1$-handle to $W'$ we obtain a smooth, compact, oriented
$4$-manifold $X$ with boundary
$\oline Y\#Y_{\ga,-1}$ and such that $X$ contains a closed surface of
genus~$g$ and
self-intersection $-1$ representing a generator of $H_2(X;\z)/\torsion$.
Now apply Theorem~\ref{hthm} to $X$ with both orientations and the Corollary
follows, since $h$ is additive under connected sums.

If $\ga$ is the $(p,q)$ torus knot in $S^3$, where $p,q$ are mutually prime
integers $\ge2$, then $Y_{\ga,-1}$ is diffeomorphic to the Brieskorn sphere
$\Si(p,q,pq-1)$. In this case one can also apply
Theorem~\ref{hthm} to the minimal
resolution of the corresponding Brieskorn singularity to get a lower bound
on $h(Y_{\ga,-1})$. In general, this lower bound does not coincide with
the upper bound given by Corollary~\ref{genus-bound}. However, we will show
in the next section that they do coincide when $p=2$, allowing us
to compute
\[h(\Si(2,2k-1,4k-3))=\lfloor k/2\rfloor\]
for $k\ge2$, where $\lfloor x\rfloor$ is the largest integer $\le x$.

\tableofcontents

\section{The general inequality}\label{main-thm}

We will now state the main result of this paper, which is more
general than Theorem~\ref{hthm}.
We first recall some definitions from \cite{Fr3}.

By a {\em lattice} we shall mean a finitely generated free abelian group
$\cl$ with a non-degenerate symmetric bilinear form $\cl\times \cl\to\z$.
The dual lattice $\text{Hom}(\cl,\z)$ will be denoted $\cl^\#$.

\be{defn}\label{eta_def}
Let $\cl$ be a (positive or negative) definite lattice. A vector
$w\in \cl$ is called {\em extremal} if $|w^2|\le|z^2|$ for all $z\in w+2\cl$.
If $w\in \cl$, $a\in \cl^\#$, and $m$ is a non-negative integer
satisfying $w^2\equiv m\mod2$ set
\[\eta(\cl,w,a,m)=\sum_z(-1)^{((z+w)/2)^2}(a\cdot z)^m,\]
where the sum is taken over all $z\in w+2\cl$
such that $z^2=w^2$.
If $m=0$ then we interpret $(a\cdot z)^m=1$ and
write $\eta(\cl,w)=\eta(\cl,w,a,m)$.
\end{defn}
Notice that $(-1)^{((z+w)/2)^2}(a\cdot z)^m$ remains the same when $z$
is replaced with $-z$. Therefore, if $w\neq0$ then this definition of
$\eta$ differs from that in \cite{Fr3} by a factor of $2$.

\be{thm}\label{h_ex}
Let $X$ be a smooth, compact, oriented $4$-manifold with a
homology sphere $Y$ as boundary and with $b^+_2(X)\ge1$.
For $i=1,\dots,b^+_2(X)$ let
$\Si_i\subset X$ be a closed surface of self-intersection $1$,
such that $\Si_i\cap\Si_j=\emptyset$ for $i\neq j$.
Set $g=\text{genus}(\Si_1)$ and suppose $\text{genus}(\Si_i)=1$ for $i\ge2$.
Let $\cl\subset H^2(X;\z)/\torsion$ be the
sublattice consisting of all vectors which
vanish on all the classes $[\Si_i]$.
Thus, $\cl$ is a unimodular, negative definite form. Let $w\in\cl$ be an
extremal vector, let
$a\in H_2(X;\z)$, and let $m$ be a non-negative integer such that
$w^2\equiv m\mod2$ and $\eta(\cl,w,a,m)\neq0$. Then
\[h(Y)+\lceil g/2\rceil+b^+_2(X)-1\ge(|w^2|-m)/4.\]
\end{thm}

The author does not see any natural generalization of the theorem in which
the roles of the $\Si_i$ are symmetric. (See Section~\ref{ns-ends} for
an explanation of this.)

\be{defn}For any definite lattice $\cl$ let $\exot(\cl)$ be the
supremum of the set of all integers $\lceil(|w^2|-m)/4\rceil$
where $w$ is any extremal vector in $\cl$ and $m$ any non-negative integer
such that (i) $w^2\equiv m\mod2$, and (ii) $\eta(\cl,w,a,m)\neq0$
for some $a\in\cl^{\#}$.
\end{defn}

With this definition of $\exot(\cl)$, the inequality in Theorem~\ref{hthm}
follows from Theorem~\ref{h_ex}. Since every unimodular lattice can be
realized as the intersection form of a smooth, compact, oriented $4$-manifold
with a homology sphere as boundary, Theorem~\ref{h_ex} implies that
$\exot(\cl)$ is finite for any unimodular definite lattice $\cl$.

Applying Theorem~\ref{h_ex} to $(\cp2\#k\cpb)\setminus\text{(open ball)}$
we see that $\exot(\cl)=0$ when $\cl$ is diagonal. (Of course, one can also
verify this directly.) If $\cl$ is not diagonal then $\exot(\cl)>0$ by 
the proof of \cite[Corollary~2]{Fr3}. In a similar fashion one can show that
$\exot(-kE_8\oplus \cl)=k$ if $\cl$ is diagonal
(cfr.\ the proof of Proposition~\ref{redhn} below).

Before embarking on the proof of Theorem~\ref{h_ex} we apply it to compute
$h$ in some examples.

\be{prop}
$h(\Si(2,2k-1,4k-3))=\lfloor k/2\rfloor$ for $k\ge2$.
\end{prop}

\proof Set $Y_k=\Si(2,2k-1,4k-3)$. The minimal resolution of the corresponding
Brieskorn singularity has intersection form isomorphic to
\be{equation}\label{ga-def}
\Ga_{4k}=\{\textstyle{\sum}x_ie_i\in\R^{4k}\st\tsum x_i\in2\z;\,
2x_i\in\z;\,x_i-x_j\in\z\},
\end{equation}
where $\{e_i\}$ is an orthonormal basis for $\R^{4k}$.
This is an
even form precisely when $k$ is even. Set $\ell=\lfloor k/2\rfloor$.
It is easy to see that
$w=\sum_{i=1}^{4\ell}e_i$ is an extremal vector in $\Ga_{4k}$ with
$\eta(\Ga_{4k},w)\neq0$. Hence
$h(Y_{k})\ge|w^2|/4=\ell$ by Theorem~\ref{h_ex}.
On the other hand, $Y_k$ is also $-1$ surgery on the $(2,2k-1)$ torus knot,
which has genus $k-1$. Since $\lceil(k-1)/2\rceil=\lfloor k/2\rfloor=\ell$,
Corollary~\ref{genus-bound} gives $h(Y_k)\le\ell$.\square

The first step in our proof of Theorem~\ref{h_ex} is the following
proposition, which uses the additivity of $h$ to reduce the theorem to the
case $h(Y)=0$, at the expense of perhaps increasing $b_2^+(X)$.
In the simplest case $b_2^+(X)=1$ it is natural to ask for an
alternative proof which only uses $4$-manifolds with $b_2^+=1$.
Unfortunately, such a proof would seem to require an extensive discussion
of bubbles, which we will not take up here.

\be{prop}\label{redhn}
If Theorem~\ref{h_ex} holds when $h(Y)=0$ and $b_1(X)=0$ then it holds
in general.
\end{prop}
\proof Performing surgery on a set of loops in $X$ representing a basis for
the free part of $H_1(X;\z)$ yields a $4$-manifold with 
$b_1=0$ but with the same intersection form and the same torsion in
$H_1(\,\cdot\,;\z)$. Since the loops can be chosen
disjoint from the surfaces $\Si_i$, it suffices to prove the Theorem when
$b_1(X)=0$.

Let $X_1$ denote the
(negative definite) $E_8$-manifold with boundary $S=\Si(2,3,5)$. Let 
$X_2$ denote the oriented $4$-manifold described by the right-handed trefoil
in $S^3$ with framing $+1$. Then $\prtl X_2=\oline S$ and $X_2$ contains
an embedded torus of self-intersection $+1$ which represents
a generator of $H_2(X_2;\z)=\z$.

For any pair $Z_1,Z_2$ of oriented
$4$-manifolds with non-empty and connected boundaries let $Z_1\#_\partial Z_2$
denote their ``boundary sum'', formed by joining the boundaries
of $Z_1$ and $Z_2$ by a $1$-handle (respecting the orientations of $Z_1$ and
$Z_2$).

If $h(Y)=-k<0$ we replace $X$ in the Theorem by the $k$-fold boundary connected
sum $X'=X\#_\prtl kX_1$. Then $\prtl X'=Y\#kS$ has $h=0$.
Also, replace $w$ by $w'=w+\sum_{j=1}^kq_j$, where $q_j$ is supported on the
$j$'th copy of $X_1$ and is given by $q_j=e_1+e_2+e_3+e_4\in E_8=\Ga_8$, in the
notation of \Ref{ga-def}. Since $\eta(E_8,e_1+e_2+e_3+e_4)=16$,
we have altered $\eta(\cl,w,a,m)$ by a factor of
$16^k$. This reduces the Theorem to the case $h(Y)\ge0$.

If $h(Y)=k>0$ replace $X$ by $X'=X\#_\prtl kX_2$, so
$\prtl X'=Y\#k\oline S$,
which has $h=0$. The $k$ copies of $X_2$ provide $k$ embedded tori in $X'$
with self-intersection $1$.\square

\section{Outline of proof}

It remains to prove Theorem~\ref{h_ex} in the case $h(Y)=0=b_1(X)$. In this
outline we will assume $Y=S^3$, since the proof in the case $h(Y)=0$ is 
almost the same. We can then just as well assume $X$ is closed.
We will also take $b^+_2(X)=1$, $H_1(X;\z)=0$, and $m=0$.
What we need to prove then is that 
if $w\in\cl$ is an extremal vector such that $w^2$ is even and
$\eta(\cl,w)\neq0$ then
\[2q\ge n+1,\]
where $q=\lceil g/2\rceil$ and $n+1=|w^2|/2$.

Let $N$ be a closed tubular 
neighbourhood of $\Si$ and $W$ the result of replacing 
$N\subset X$ by $V=N\#\cpb$. 
After choosing orientations of $\Si$ and $S=\cp1\subset\cp2$ we form
two internal connected sums $\Si^\pm\subset V$ of $\Si$ and $S$, in one case
preserving orientations and in the other case reversing them. 
Thus $\Si^\pm$ has self-intersection number $0$. We
construct a smooth $1$-parameter family of metrics $g(t)$ on $W$, independent
of $t$ outside $V$, as follows:
First choose an initial metric in which a small tubular
neighbourhood of $C=\prtl V$ is
isometric to $[0,T]\times C$, where $T\gg0$.
Then stretch this initial metric along the
link of $\Si^+$ (for $t>1$) and along the link of $\Si^-$ (for $t<-1$).
More precisely, when $t>1$ then $W$ should contain a cylinder
$[0,t]\times S^1\times\Si^+$, and similarly for $t<-1$. 

Let $s\in H^2(W;\z)$ be the Poincar\'e dual of $[S]$. Choose an integral
lift $w'$ of $w$ and for any integer $k$
let $E_k\to W$ be the $\U2$ bundle with $c_1=w'+s$ and
$c_2=k$. For each $t$ let $M_{k,t}$ be the
moduli space of projectively $g(t)$ anti-self-dual connections in $E_k$. Let
$\cm_k$ denote the disjoint union of the $M_{k,t}$, $t\in\R$.
For suitable, small
perturbations of the ASD equations, the irreducible part $\cm^*_k$ of $\cm_k$
will be a smooth manifold of dimension 
\[\dim\,\cm^*_k=8k+4n+1.\]
In Section~\ref{sec:red} we will show that $\cm_k$ contains no
reducibles for $k<0$, while $\cm_0$ contains a finite number of reducibles
(for a generic choice of initial metric).
Moreover, one can describe explicitly, in terms of
$\cl$ and $w$ only, the splittings of $E_0$ (into complex line bundles)
that correspond to reducibles in $\cm_0$. To each such
splitting there is associated a ``degree'' (always $\pm1$)
which measures the number of times (counted with sign) that this
splitting occurs. In addition to the degree, each splitting also comes with
a sign which measures whether the overall orientation of a certain determinant
line bundle over the orbit space of connections in $E_0$ agrees with the
``complex orientation'' at the corresponding reducible points.
The sum of (degree)$\cdot$(sign) over all the splittings
is equal to $2\eta(\cl,w)$.

Let $\cm_0'\subset\cm_0$ be the result of removing from $\cm_0$
a small open neighbourhood about each reducible point. For $k<0$ set
$\cm_k'=\cm_k$.
Recall that to any base-point in $W$ one can associate
a principal $\SO3$ bundle $\e\to\cm'_k$ called the base-point fibration.
For $j\ge0$ and any subset $S$ of $\cm'_k$
let $x^j\cdot S$ denote the intersection of $S$ with $j$
generic geometric representatives for $-\frac14p_1(\e)$, in the sense of 
\cite{KM3}. These representatives should depend only on the restriction of
elements of $S$ to suitable compact subsets of $W\setminus V$ where
the metric does not vary with $t$.

We also need a variant of this construction
where the location of the base-point depends on $t$; in this case the
result is denoted $x_1^j\cdot S$. To make this more precise: for
$\pm t\gg0$ the geometric representatives should be defined through restriction
to subsets of some tubular neighbourhood of $\Si^\pm$ where the metric does
not vary, while for intermediate values of $t$ we interpolate,
in a certain sense.

Now suppose the theorem does not hold, ie that $2q\le n$, and set
\[\hcm=\text{$1$-dimensional part of $\;x^{n-2q}\cdot(x_1^2-4)^q\cdot\sum_{k\le0}
\cm'_k$.}\]
This is a formal linear combination of oriented $1$-manifolds with boundary.
The number of boundary points of $\hcm$, counted with multiplicity, equals
$\pm4^{-n}\eta(\cl,w)$. On the other hand, it follows from Mu\~noz'
description of the ring structure of the Floer cohomology
of $S^1\times\Si^\pm$ that
the number of ends of $\hcm$, counted with multiplicity, is zero.
This contradiction proves the theorem in the case considered.

\section{Parametrized moduli spaces}\label{par-mod}

In this section we study instanton moduli spaces over a $4$-manifold $W$ 
with tubular ends, parametrized by a family of tubular end metrics.
We are interested in questions of orientations and reducibles.

\subsection{Orientations}

Let $W$ be an oriented Riemannian $4$-manifold with tubular ends
$\rpy_j$, $j=1,\dots,r$, so $W\setminus\cup_j(\rpy_j)$ is compact.
Let $E\to W$ be a $\U2$ bundle. We denote by $\frakg_E$
the bundle of Lie algebras
associated to $E$ and the adjoint representation of $\U2$ on its Lie algebra
$\u2$, and by $\frakg'_E$ the sub-bundle corresponding to the sub-algebra
$\su2\subset\u2$. 

Choose an isomorphism $E|_{\rpy_j}=\R\times E_j$, where $E_j$ is a $\U2$
bundle over $Y_j$. Let $\al_j$ be a connection in 
$E_j$ such that the induced connection in $\frakg'_{E_j}$
is non-degenerate flat. (As usual, `flat connections' will in practice mean
critical points of the perturbed Chern-Simons functional as in \cite{F1}.)
Let $\si_j\ge0$ be small, and $\si_j>0$ if
$\al_j$ is reducible. Fix an even integer $p>4$.
Choose a smooth connection $A_0$ in $E$
which agrees with the pull-back of $\al_j$ over the $j$'th end, and
define the Sobolev space of connections
\[\ca=\ca(E,\al)=\{A_0+a\st a\in\llw p\si1(T^*W\otimes\frakg'_E)\}.\]
Here $\llw p\si k$ is the space of sections $s$ such that
$(\nabla_{A_0})^j(e^ws)\in L^p$ weakly for $0\le j\le k$, where 
$w:W\to\R$ is any smooth function with $w(t,y)=e^{\si_jt}$ for
$(t,y)\in\rpy_j$.
As explained in \cite{D5} there is a Banach Lie group $\cg$ consisting of
$L^p_{2,\loc}$ gauge transformations, such that $\cg$ acts smoothly on $\ca$ 
and the following holds: 
If $A,B\in\ca$ and there is an
$L^p_{2,\loc}$ gauge transformation $u$
such that $u(A)=B$ then $u\in\cg$.


There is a real determinant line bundle
$\det(\del)$ over $\cb=\ca/\cg$ associated to the family of Fredholm operators
\[\del_A=d^*_A+d^+_A:\lls p1\to\ls p,\]
see \cite{DK,D5} and Appendix~\ref{sec:detlb}.
It is proved in \cite{D5} that $\det(\del)$ is orientable. 
An orientation of $\det(\del)$ defines an orientation 
of the regular part of the instanton moduli space
$M^*\subset\cb^*$ cut out by the equation $F^+_0(A)=0$, where $F_0(A)$
is the curvature of the connection that $A$ induces in $\frakg'_E$.
As usual, $\cb^*$ denotes the irreducible part of $\cb$, etc.
As for the choice of orientation, the important thing for us will be that
orientations of the various moduli spaces involved be chosen compatible
with gluing maps, and this can be done at least in the situations we will
consider.

More generally, let $g(t)$ be a smooth family of Riemannian metrics on
$W$, where $t$ runs through some parameter space $\R^b$, such that
$g(t)$ is constant in $t$ outside some compact set in $W$.
In the following let 
$\La^j$ be the bundle of $j$-forms and $\La^+_t$ the bundle of
$g(t)$ self-dual $2$-forms over $W$.
There is then a parametrized moduli space
\[\cm\subset\cb\times\R^b\]
consisting of all pairs $([A],t)$ satisfying
\be{equation}\label{asdpar}
P_\bt^+F_0(A)=0,
\end{equation}
where $P_t:\La^2\to\La^+_t$. 
Generically, the irreducible part of this moduli space, $\cm^*$,
will be a smooth manifold. (In the situations encountered in this paper
it will suffice to
start with the usual translationary invariant perturbations over
the ends of $W$, given by perturbations of the Chern-Simons functional
as in \cite{F1}, and add further 
perturbations defined in terms of holonomy along a finite number
of thickened loops in $W$. In general, perturbations should be small so that
one can control reducibles.)

To orient $\cm^*$ we note that
if $\bt$ is sufficiently close to a reference point $\tau$ then 
\Ref{asdpar} is equivalent to
\[P_\tau^+P_\bt^+F_0(A)=0.\]
The derivative of $(A,t)\mapsto P^+_\tau P^+_tF_0(A)$ at $(A,\tau)$ is
the operator
\be{align*}
T_{A,\tau}:\lls p1(\La^1\otimes\frakg'_E)\times\R^b&\to
\ls p(\La^+_\tau\otimes\frakg'_E)\\
(a,x)&\mapsto P^+_\tau d_Aa + S_{A,\tau}x,
\end{align*}
where for any $(A,\tau)\in\ca\times\R^b$ we define
\[S_{A,\tau}x=\left.\frac\prtl{\prtl s}\right|_0P^+_\tau
P^+_{\tau+sx}F_0(A).\]
Differentiating the equation $P^+_\tau P^-_\tau=0$ we obtain
\[\frac\prtl{\prtl s}
(P^+_{\tau+sx}P^-_\tau-P^+_\tau P^+_{\tau+sx})=0,\]
hence
\[S_{A,\tau}x=\left.\frac\prtl{\prtl s}\right|_0P^+_{\tau+sx}F_0(A)\qquad
\text{if $([A],\tau)\in\cm$.}\]
A point $([A],t)\in\cm$ is called {\em regular} if $T_{A,t}$ is surjective.
The tangent space of $\cm$ at an irreducible regular point $([A],\bt)$ can
be identified with the kernel of the operator
\[\ti\del_{A,\bt}=\del_{A,t}\oplus S_{A,t}:
\lls p1(\La^1\otimes\frakg'_E)\oplus\R^b\to
\ls p(\frakg'_E)\oplus\ls p(\La^+_t\otimes\frakg'_E),\]
where $\del_{A,t}$ denotes the operator $\del_A$ in the metric $g(t)$.
What we are seeking, therefore, is an orientation of the determinant line
bundle $\det(\ti\del)$ over $\cb\times\R^b$ associated to
the family of Fredholm operators $\ti\del_{A,\bt}$.
But as explained in Subsection~\ref{cont-det} there is a canonical
isomorphism
\be{equation}\label{detdel-iso}
\ga:\dett(\del)\to\dett(\ti\del))\otimes(\dett(\R^b))^*.
\end{equation}
Given an orientation of $\det(\del)$,
this orients $\det(\ti\del)$, hence $\cm^*$.

\subsection{Local structure near reducibles}
\label{loc-str}

We now assume, for simplicity, that each $Y_j$ is a rational homology sphere,
and that $b_1(X)=0$. More importantly, we take $b=b^+$, where $b^+=b_2^+(W)$.

Let $D=([A],\bt)\in\cm$ be a reducible point, so $A$ respects some splitting
of $E$ into complex line
bundles, $E=L_1\oplus L_2$.
Then $\frakg_E'=\uline{\R}\oplus K$, where $K=L_1\otimes\bar L_2$
and $\uline{\R}=W\times\R$. Here the constant section $\mathbf1$ of $\uline\R$
acts as $2i$ on $K$.
Note that changing the order of $L_1,L_2$ changes the
orientation of $K$, hence many signs in what follows.

Let $B$ denote the connection that $A$ induces in $K$. Noting that
$F_0(A)$ takes values in $\uline\R$, we can identify $F_0(A)=-2iF(B)$,
which is an anti-self-dual, closed $L^2$ $2$-form representing $4\pi c_1(K)$.

Let $\harm^+_\bt$ be the space of self-dual, closed, $L^2$ $2$-forms on
$(W,g(t))$. By \cite[Proposition~4.9]{APS1}
we can identify $\harm^+_t$ with a subspace of $H^2(W)$. 

The operator
\[P^+_td_A:\lls p1(\La^1\otimes \frakg'_E)\to\ls p(\La^+_t\otimes \frakg'_E)\]
is the sum of two operators: one with values in $\uline\R$, which is the usual
$P^+_td$ operator, and another with values in $K$, which we call $P^+_td_B$.
Note that the dual of the cokernel of $P^+_td$ is $\harm^+_\bt$,
and the map $\R^{b^+}\to\coker(\del_{A,t})$ defined by $S_{A,t}$
takes values in $(\harm^+_t)^*$.
Let $R_D:\rbp\to(\harm^+_t)^*$ denote this operator.
\be{observation}
$D$ is a regular point of $\cm$ if and only if the following two
conditions hold:
\be{description}
\item[(i)]$P^+_td_B:\lls p1(\La^1\otimes K)\to\ls p(\La^+_t\otimes K)$
is surjective.
\item[(ii)]$R_D:\rbp\to(\harm^+_t)^*$ is an isomorphism.
\end{description}
\end{observation}

Note that if the bundle $K$ is non-trivial then, in the metric $g(t)$,
$d^+_B$ is surjective if and only if $d^*_B+d^+_B$ is surjective,
in which case the latter operator has non-negative index. So if this index
is negative then $D$ cannot be a regular point of $\cm$. On the other hand,
if the index is non-negative (which will be the case in our applications)
then as explained in \cite{D1} there is a simple
local perturbation of the anti-self-duality equation near $D\in\cb\times\rbp$
such that $D$ solves the perturbed equation and such that the perturbed 
analogue of $P^+_td_B$ is surjective.
(Ideally, one should look for a generalization of Freed and Uhlenbeck's theorem
to our situation, but we will not pursue this here.)

We will now give a more concrete description of the map $R_D$. 
Let $\Si_1,\dots,\Si_{b^+}\subset W$ be closed surfaces such that
$\Si_i\cdot\Si_j=0$ if $i\neq j$ and $\Si_i\cdot\Si_i>0$ for all $i$.
Then there is a linear ismorphism
\[\harm^+_t\to\R^{b^+},\quad\om\mapsto
\left(\int_{\Si_i}\om\right)_{1\le i\le b^+}.\]
We can therefore define a basis $\{\om_{i,t}\}\subset\harm^+_t$
by the conditions
\[\int_W\om_{i,t}\wedge\om_{i,t}=1;\quad\int_{\Si_i}\om_{i,t}>0;\quad
\int_{\Si_i}\om_{j,t}=0
\quad\text{for $i\neq j$.}\]
For any $v\in H^2(W)$ consider the map
\be{equation}\label{fvdef}
f_v:\R^{b^+}\to\R^{b^+},\quad\bt\mapsto
\left(\int_W\om_{i,t}\wedge v\right)_{1\le i\le b^+},
\end{equation}
where in the integral $v$ is represented by some bounded, closed $2$-form.
Clearly, $f_v(t)=0$ if and only if $v$ can be represented by a $g(t)$
anti-self-dual, closed $L^2$ form. 

\be{prop}With respect to the basis for $\harm^+_t$ defined above we have
\[R_D=df_v(t),\]
where $v=4\pi c_1(K)$, and $df_v(t)$ is the derivative of the
function $f_v$ at $t$.
\end{prop}
\proof For every $x\in\rbp$ we have
\be{align*}
R_Dx\cdot\om_{i,\bt}
&=\left.\frac d{ds}\right|_0\int_WP^+_{\bt+sx}F_0(A)\wedge\om_{i,\bt}\\
&=\left.\frac d{ds}\right|_0\int_WP^+_{\bt+sx}F_0(A)\wedge\om_{i,\bt+sx}\\
&=\left.\frac d{ds}\right|_0\int_Wv\wedge\om_{i,\bt+sx}\\
&=df_{v,i}(t)x,
\end{align*}
where $f_{v,i}$ is the $i$'th component of $f_v$.
\square

We conclude this section with a simple result about orientations. Suppose
$D$ is a regular point. Then $\ker(\del_D)=\ker(\ti\del_D)$ is a 
complex vector space (recall that we assume $b_1(W)=0$),
and therefore has a canonical orientation. Furthermore,
\[\coker(\ti\del_D)=\R\mathbf1,\quad
\coker(\del_D)=\R\mathbf1\oplus(\harm^+_t)^*.\]
Given the ordering of the surfaces $\Si_i$ we then obtain orientations of
$\det(\del_D)$ and $\det(\ti\del_D)$, which we refer to as the
``complex orientations'', cfr.\ \cite{D2}.
We would like to know how these compare under the isomorphism $\ga_D$ in
\Ref{detdel-iso}. As explained in Subsection~\ref{subsec:detl} the 
exact sequence
\be{equation*}
0\to\ker(\del_D)\to\ker(\ti\del_D)\to\R^b\to
\coker(\del_D)\to\coker(\ti\del_D)\to0
\end{equation*}
gives rise to a natural isomorphism
\[\ga'_D:\det(\del_D)\to\det(\ti\del_D)\otimes(\det(\R^b))^*.\]
By Proposition~\ref{gagap} (with $\eps_1=\eps_2=1$) we have
\be{equation}\label{gadd}
\ga_D=\ga'_D.
\end{equation}
This yields:
\be{prop}\label{algebra}
If $D$ is a regular point then the complex orientations of
$\det(\del_D)$ and $\det(\ti\del_D)$ agree under the isomorphism
$\ga_D$ if and only if $R_D$ preserves orientation.
\end{prop}

\section{The family of metrics}\label{fam-metrics}

We can now begin the proof of Theorem~\ref{h_ex} in earnest.
By Proposition~\ref{redhn} we may assume $b_1(X)=0$. (The assumption
$h(Y)=0$ will not be used until Section~\ref{main-mech}.) In this section
we construct from $X$ the $4$-manifold $W$ that will be the base-manifold
for our moduli spaces later.
We obtain a $b_2^+(W)+1$ dimensional family of Riemannian
metrics on $W$ by stretching
along various hypersurfaces.

To define $W$, set $b^+=b_2^+(X)$ and let $\hat X$ be the result of adding
a half-infinite tube $\R_+\times Y$ to X.
Choose disjoint, compact tubular neighbourhoods $N_i$ of the surfaces $\Si_i$
and let $V_i\approx N_i\#\cpb$ be the blow-up of $N_i$ at some interior
point away from
$\Si_i$. We then define $W$ to be the manifold obtained from
$\hat X$ by replacing each $N_i$ by $V_i$. Thus,
\[W\approx\hat X\#b^+\cpb.\]
Let $S_i\subset\cp2$ be a sphere representing a generator of $H_2(\cp2;\z)$.
Choose internal connected sums $\Si^+_i=\Si_i\#S_i$ and
$\Si^-_i=\Si_i\#\oline S_i$ in the interior of $V_i$.

Let $S^1$ be the boundary of the closed unit disk $D^2\subset\R^2$ centered
at the origin. For each $i$ choose smooth embeddings
\[q^\pm_i:D^2\times\Si^\pm_i\to V_i\]
such that $q^\pm_i(0,\,\cdot\,)$ is the identity on $\Si^\pm_i$.
Let $N^\pm_i$ denote the image of the embedding
\be{align*}
\rho^\pm_i:&[0,1]\times S^1\times\Si^\pm_i\to V_i\\
&(s,x,z)\mapsto q^\pm_i(\frac12(s+1)x,z).
\end{align*}

Set $W'=W\setminus(\cup_i\,\text{int}\,V_i)$ and $C=\prtl W'$, and 
let $N_C$ be the image of a smooth embedding
\[\rho_C:[0,1]\times C\to W'\]
such that $\rho_C(0,\,\cdot\,)$ is the identity on $C$. Set
$W^-=W'\setminus\rho_C([0,1)\times C)$.
Choose Riemannian metrics on $C$ and on
$S^1\times\Si^\pm_i$.

Choose a smooth function $\tau:\R\to\R$ such that
$\tau'\ge0$, $\tau(s)=0$ for $s\le\frac13$, $\tau(s)=1$ for $s\ge\frac23$.
For $r\ge1$ let $\xi_r:[0,r]\to[0,1]$ be the diffeomorphism whose inverse
is given by
\[\xi\inv_r(s)=s+(r-1)\tau(s).\]

Now choose a smooth family of Riemannian metrics $g(T,t)$ on $W$, where
$T\ge1$ and $t=(t_1,\dots,t_{b^+})\in\rbp$, such that the following holds:
\be{itemize}
\item If $N_C$ and $N^\pm_i$ have the metrics induced by $g(T,t)$,
intervals in $\R$ have the standard metric, and
$S^1\times\Si^\pm_i$ and $C$ have the metrics chosen above
(which are independent of $T,t$),
then the following composite maps should be isometries:
\be{align*}
&[0,\pm t_i]\times S^1\times\Si^\pm_i\oset{\xi_{\pm t_i}\times\text{Id}}
\longrightarrow[0,1]\times S^1\times\Si^\pm_i\oset{\rho^\pm_i}\longrightarrow
N^\pm_i\qquad\text{if $\pm t_i\ge1$}\\
&[0,T]\times C\oset{\xi_T\times\text{Id}}\longrightarrow
[0,1]\times C\oset{\rho_C}\longrightarrow N_C.
\end{align*}
\item $g(T,t)$ is independent of $T$ outside $N_C$.
\item $g(T,t)$ is independent of $t_i$ outside $N^+_i\cup N^-_i$.
\item For $t_i\ge1$, $g(T,t)$ is independent of $t_i$ outside $N^+_i$.
\item For $t_i\le-1$, $g(T,t)$ is independent of $t_i$ outside $N^-_i$.
\item $g(T,t)$ is on product form on $\rpy$.
\end{itemize}

As $T\to\infty$ we obtain from $(W,g(T,t))$ the following
Riemannian manifolds with tubular ends:
\be{align*}
W^-_\infty&=W^-\cup_C(\R_-\times C)\\
V_{i,\infty}&=V_i\cup_{\prtl V_i}(\R_+\times\prtl V_i)
\end{align*}
The metric on $W^-_\infty$ is independent of $t$, while the metric
on $V_{i,\infty}$ depends only on $t_i$ and is denoted $g_i(t_i)$.

\section{Reducibles}\label{sec:red}

Let $s_i\in H^2(W)$ be the Poincar\'e dual of $[S_i]$, and 
choose a $U(2)$ bundle $E\to W$ such that $c_1(E)=w+\sum_is_i$ modulo
torsion, where $w$ is as in Theorem~\ref{h_ex}.
For any integer $k$ let $M\tt(E,k)$ denote the moduli space of
projectively $g(T,t)$ anti-self-dual connections $A$ in $E$ which are
asympotically trivial over the end $\rpy$ and has ``relative second Chern
class'' $k$, ie
\be{equation}\label{e-est}
\fact\int_W\text{tr}(F_A^2)=k-\frac12c_1(E)^2,
\end{equation}
where $F_A$ is the curvature of $A$.
We wish to determine the reducible connections in $\bigcup_\bt M\tt(E,k)$ when
$T$ is large and $k\le0$.

For any $v\in H^2(W)$ let $f_{T,v}:\rbp\to\rbp$ be the map defined in 
Subsection~\ref{loc-str}, using the family of metrics $g(T,t)$, $t\in\rbp$.
If $v\cdot[\Si^\pm_i]\neq0$
for all $i$ and both signs, 
then the zeros of $f_{T,v}$ form a bounded set. Restricting 
$f_{T,v}$ to a large sphere, of radius $r>0$ say, we then obtain a map
\[f_{T,r,v}:S^{b^+-1}_r\to\R^{b^+}\setminus0\]
and we define
\be{equation}\label{eq:degf}
\deg(v)=\deg(f_{T,r,v}/|f_{T,r,v}|).
\end{equation}
By the homotopy invariance of the degree, the right hand side is independent
of $T$ and $r\gg0$, so $\deg(v)$ is well-defined.

\be{lemma}\label{Tred}
For $T>0$ sufficiently large the following holds.
\be{description}
\item[(i)]If $k<0$ then there are no reducibles in
$M\tt(E,k)$ for any $\bt$.
\item[(ii)]Suppose $M\tt(E,0)$ contains a reducible connection
which respects a
splitting $E=L_1\oplus L_2$. Then modulo torsion the Chern class
$c_1(L_1\otimes\bar L_2)$ has the form
$v=z-\sum_i\eps_is_i$, where $z\in w+2\cl$, $z^2=w^2$,
$\eps_i=\pm1$. Moreover, for such $v$ we have
\[\deg(v)=\eps_1\eps_2\dots\eps_{b^+}.\]
\end{description}
\end{lemma}
\proof Fix $k\le 0$. If $[A]\in M\tt(E,k)$ respects the splitting
$E=L_1\oplus L_2$ and $B$ is the induced connection in $K=L_1\otimes\bar L_2$
then $\phi=\frac1{2\pi i}F_B$ is an anti-self-dual, closed $L^2$ form
(with respect to $g(T,t)$) which represents $v=c_1(K)$. Recall from
\cite[Section~4]{Fr3} that
\[c_1(E)^2=v^2+4k.\]
Now suppose $[A_n]$ is a reducible point in $M_{T(n),t(n)}(E,k)$
for $n=1,2,\dots$, where $T(n)\to\infty$ as $n\to\infty$.
Since $\|\phi_n\|_{L^2}$ is independent of $n$, by the above equation, and
$c_1(E)\cdot[\Si^\pm_i]=\mp1$ for all $i$, $\{t(n)\}$ must be a bounded
sequence. After passing to a subsequence we may therefore assume that
$\{t(n)\}$ converges in $\rbp$, and that
$\{\phi_n\}$ converges in $C^\infty$ over compact subsets of both
$W^-_\infty$ and $V_{i,\infty}$, with limits $z'$ and $z'_i$, respectively.

Note that $z',z'_i$ are harmonic $L^2$ forms and therefore decay
exponentially over the ends (see \cite{APS1}).
We identify the cohomology class $[z']$ with the element $z\in H^2(W)$
which maps to $(0,[z'])$ under the isomorphism
\be{equation}\label{c-iso}
H^2(W)\to H^2(W^+)\oplus\ker[H^2(W^-)\to H^2(C)],
\end{equation}
where $W^+=W\setminus\text{int}\, W^-\approx\cup_iV_i$. Then
$z\in w+2\cl$, so $z^2\le w^2$ since $w$ is extremal in $\cl$. 
We also identify $[z'_i]$ with the class $z_i\in H^2(W)$ which maps to
$([z'_i],0)$ in~\Ref{c-iso}. Then $z_i=\al_i\si_i+\beta_is_i$,
where $\al_i$ is an even integer and $\beta_i$ is an odd integer.
Anti-self-duality implies
$z_i^2\le0$, so $z_i^2=\al_i^2-\beta_i^2\le-1$.
We therefore get, for large $n$,
\[w^2-b^+=c_1(E)^2=v_n^2+4k\le z^2+\sum_iz_i^2+4k\le w^2-b^++4k.\]
This gives $k\ge0$, with equality only if $z^2=w^2$ and
$\al_i=0$, $\beta_i=\pm1$ for each $i$.

We will now compute the degree of $v=z-\sum_i\eps_is_i$.
Let $\harm^+_{T,t}$ be the space of self-dual, closed, $L^2$ $2$-forms
on $W$ with respect to the metric $g(T,t)$, and let $\om_{i,T,t}$,
$i=1,\dots,b^+$ be the basis for $\harm^+_{T,t}$ constructed in
Subsection~\ref{loc-str}, using the surfaces $\Si_i$.
Note that as $T\to\infty$, $\om_{i,T,\bt}$ converges in $C^\infty$
to zero over compact subsets of $W^-_\infty$ (since the intersection form
on $\ker[H^2(W^-)\to H^2(C)]$ is negative definite; see
\cite[Proposition~4.9]{APS1})
and also over compact subsets of $V_{j,\infty}$ for $j\neq i$,
and it converges over compact subsets of $V_{i,\infty}$ to some
$g_i(t_i)$ self-dual form $\eta_i=\eta_i(t_i)$ uniquely determined by the
properties
\[\int_{\Si_i}\eta_i\ge1;\quad\int_{V_{i,\infty}}\eta_i\wedge\eta_i=1.\]
Moreover, the convergence is uniform for $|\bt|\le r$, so $f_{T,v}|_{D(r)}$
converges in $C^0$ to $f_{v,1}\times\dots\times f_{v,b^+}$
as $T\to\infty$, where
\[f_{v,i}:\R\to\R,\quad s\mapsto\int_{V_{i,\infty}}\eta_i(s)\wedge v.\]
Now,
\[\int_{\Si_i}\eta_i(s)\pm\int_{S_i}\eta_i(s)=
\int_{\Si_i^{\pm}}\eta_i(s)\to0\quad\text{as $s\to\pm\infty$.}\]
Hence $\int_{S_i}\eta_i(s)$ is negative for $s\gg0$ and 
positive for $s\ll0$. It follows easily from this that
$\deg(v)=\eps_1\eps_2\dots\eps_{b^+}$.
\square

In the next section we will need to understand
which instantons over $W$ may restrict to reducible connections over
open subsets. We record here the following basic result.

We say a connection $A$ in an $\SO3$ bundle $P$ is {\em reducible} if
$A$ is preserved by a non-trivial automorphism of $P$, and we say $A$ is
{\em s-reducible} if $A$ is preserved by a 
non-trivial automorphism of $P$ which lifts to $P\times_{\text{Ad}}\SU2$,
cfr.\ \cite[Section~2]{Fr3}.
We denote by $\frakg_P$ the bundle of Lie algebras associated to $P$.

\be{lemma}[\cite{DK,KM3}]\label{twred}
Let $Z$ be a connected, oriented Riemannian $4$-manifold and $A$ an
anti-self-dual connection in a principal $\SO3$ bundle $P\to Z$.
Suppose $A$ is not flat, and that $A$ is reducible over some non-empty
open subset of $Z$. Then
$A$ respects some splitting $\frakg_P=\lla\oplus L$ where $\lla$ is a real
line bundle and $L$ a real $2$-plane bundle. Moreover, $A$ is s-reducible if 
and only if $\lla$ is trivial.
\end{lemma}

\proof Let $S$ be the set of points $x\in Z$ such that $A$ is reducible 
in some open neighbourhood of $x$. Then $S$ is open, and non-empty by
assumption. The proof of \cite[Lemma~4.3.21]{DK} shows that $S$ is also
closed, hence $S=Z$. Therefore, $A$ locally respects splittings of the form
$\frakg_P=\lla\oplus L$, where $F_A$ takes values in $\lla$.
Since $(d^*_Ad_A+d_Ad^*_A)F_A=0$, unique continuation (see \cite{JKazdan1})
implies that
$F_A$ cannot vanish in any non-empty open subset of $Z$. It follows that
$A$ respects a global splitting.
The last assertion of the lemma is left to the reader.\square

\section{Cutting down parametrized moduli spaces}

From now on we fix a large $T>0$ such that the conclusions of Lemma~\ref{Tred}
hold, and suppress $T$ from notation.

For any projectively flat connection $\rho$ in $E|_Y$ set
\[\cm_\rho=\bigcup_{t\in\rbp}(M_t(E,\rho)\times\{t\}),\]
where $M_t(E,\rho)$ is the moduli space of projectively $g(t)$ anti-self-dual
connections in $E$ with limit $\rho$ over the end $\rpy$.
We will now explain what we shall mean by ``cutting down'' $\cm^*_\rho$
(or more generally, a subset of $\cm^*_\rho$) according to a monomial
$\kappa\prod_ix^{n_i}_i$, where $\kappa=z_1\cdots z_J$, $z_j\in H_{d_j}(W^-;\z)$,
$0\le d_j\le2$, and each $n_i$ is a non-negative integer. Each
$x_i$ may be thought of as the point class in $H_0(W;\z)$, but 
the location of the point will depend on $t_i$.
The cut down moduli space will be denoted
\be{equation}\label{cdms}
(\kappa\tprod_ix_i^{n_i})\cdot\cm^*_\rho,
\end{equation}
although it depends on various choices not reflected in the notation.

Roughly speaking, to cut down $\cm^*_\rho$ according to
$\kappa$ we first restrict connections to $W^-$
(where the metric does not vary), and then cut down by the product of
the $\mu$-classes of the $z_i$ just as on a closed $4$-manifold.
As for the factor
$x^{n_i}_i$, let $\cm^*_{\pm t_i\ge r}$ be the part of $\cm^*_\rho$ 
where $\pm t_i\ge r$, and let $r\gg0$.
To cut down $\cm^*_{\pm t_i\ge r}$ according to
$x^{n_i}_i$ we restrict connections to a tubular neighbourhood of 
$\Si^\pm_i$ where the metric does not vary, and then cut down by the
$n_i$'th power of the $\mu$-class of a point. In the intermediate region
$|t_i|\le r$ we interpolate by ``moving base-points'', as we will explain
in a moment.

We will now make this precise.
To cut down by $\kappa$, choose disjoint, compact, codimension~$0$
submanifolds $U_j\subset W^-$ such that $z_j$ is the image of a class
$z_j'\in H_{d_j}(U_j)$. To ensure that irreducible instantons
over $W$ restrict to irreducible connections over $U_j$ we also require that
$H_1(U_j;\z/2)\to H_1(W;\z/2)$ be surjective.
Let $\cb(U_j)$ be the orbit space of $L^p_1$ connections in $E|_{U_j}$ with
a fixed central part, and $\cb^*\subset\cb$ the irreducible part. 
Choose a generic geometric
representative $R_j\subset\cb^*(U_j)$ for $\mu(z'_j)\in H^{4-d_j}(\cb^*(U_j))$,
(see \cite{DK,KM3}). Set
\[Z_\kappa=\{([A],t)\in\cm^*_\rho\st\text{$[A|_{U_j}]\in R_j$ for
$j=1,\dots,J$}\}.\]

We now turn to the factor $x_i^{n_i}$.
For each $i$ choose disjoint closed subintervals
$\{I_{i\nu}\}_{1\le\nu\le n_i}$ of $[\frac23,1]$, each with non-empty
interior, and set
\[B^\pm_{i\nu}=\rho^\pm_i(I_{i\nu}\times S^1\times\Si^\pm_i),\]
where $\rho^\pm_i$ is as in Section~\ref{fam-metrics}. Choose disjoint
compact, connected codimension~$0$ submanifolds
$K\iu\subset W\setminus\cup_{j\neq i}V_j$ such that
\be{itemize}
\item $B^\pm\iu\subset K\iu$,
\item The sets $K\iu$ are mutually disjoint and also disjoint from the sets
$U_j$,
\item $H_1(K\iu;\z/2)\to H_1(W;\z/2)$ is surjective.
\end{itemize}
When the perturbations of the ASD equations are sufficiently small
then the restriction map
\[r\iu:\cm^*_\rho\to\cb^*(K\iu)\]
is well-defined, as follows from Lemma~\ref{twred} and a compactness argument.
Also, when $r\gg0$ there is a well-defined restriction map
\[r^\pm\iu:(\cm^*_\rho)_{\pm t_i>r}\to\cb^*(B^\pm\iu).\]
This follows by a compactness and unique continuation argument from the fact 
that there are no reducible projectively flat connections
in $E|_{S^1\times\Si^\pm_i}$ (since $c_1(E)\cdot[\Si^\pm_i]=\mp1$).

Choose a base-point $b^\pm\iu\in B^\pm\iu$ and a smooth path
$\ga\iu:[0,1]\to K\iu$ such that $\ga\iu(\pm1)=b^\pm\iu$.
Let $\e\iu\to\cb^*(K\iu)$ and $\e^\pm\iu\to\cb^*(B^\pm\iu)$ be the 
complexifications of the natural
$3$-plane bundles associtated to the base-points $b^+_i$ and $b^\pm_i$,
respectively. By means of holonomy along $\ga\iu$ we can identify $\e\iu$ with
the corresponding bundle with base-point $b^-\iu$.

For each pair $i,\nu$ choose a generic pair of
sections of $\e\iu$; pulling these back by $r\iu$ gives a pair
$(s^0_{i\nu1},s^0_{i\nu2})$ of sections of the bundle
$\f\iu=r^*\iu(\e\iu)$ over $\cm^*_\rho$. Choose also a generic pair of sections
of $\e^\pm\iu$; this gives a pair $(s^\pm_{i\nu1},s^\pm_{i\nu2})$ of
sections of $\f^\pm\iu=(r^\pm\iu)^*\e^\pm\iu$. Let $\beta:\R\to\R$ be a smooth
function such that $\beta(s)=0$ for $|s|\le r+1$ and $\beta(s)=1$ for
$|s|\ge r+2$. Because of the canonical identification $\f\iu=\f^\pm\iu$ over
$(\cm^*_\rho)_{\pm t_i>r}$, it makes sense to define sections
$(s_{i\nu1},s_{i\nu2})$ of $\f\iu$ by
\[s_{i\nu j}=
\beta(t_i)s^\pm_{i\nu j}+(1-\beta(t_i))s^0_{i\nu j}
\quad\text{for $\pm t_i\ge0$}.\]
Now let
\[Z\iu\subset\cm^*_\rho\]
be the locus where $s_{i\nu1},s_{i\nu2}$ are linearly dependent. As
explained in \cite{KM3}, $Z\iu$ is a disjoint union of finitely many
smooth submanifolds, such that the top stratum has codimension~$4$ and
a natural orientation, and there is no stratum of codimension~$5$.
We now define
\[(\kappa\tprod_ix_i^{n_i})\cdot\cm^*_\rho=Z_\kappa\cap(\cap\iu Z\iu).\]

\section{Main mechanism of proof}\label{main-mech}

We will use the following characterization of the invariant $h$, see
\cite{Fr3} for more details. For an oriented integral homology
$3$-sphere $V$ let $\cf^*(V)$ be the $\z/8$ graded Floer cochain
complex of $V$ with rational coefficients. Let $\del'\in\cf^1(V)$ and
the homomorphism $\del:\cf^4(V)\to\q$ be defined by counting with sign
the points in $0$-dimensional moduli spaces $\check M(\al,\theta)$ and 
$\check M(\theta,\al)$ over $\R\times V$,
respectively, where $\al$ is an irreducible flat
$\SU2$ connection and $\theta$ the trivial $\SU2$ connection. 
Here $\check M=M/\R$. We denote
by $\del'_0\in\hf^1(V)$ and $\del_0:\hf^4(V)\to\q$ the corresponding
data in cohomology. Then
\be{itemize}
\item Either $\del_0=0$ or $\del'_0=0$
\item $h(V)>0$ if and only if $\del_0\neq0$
\item $h(V)<0$ if and only if $\del'_0\neq0$.
\end{itemize}
By Proposition~\ref{redhn} we may assume $h(Y)=0$.
Taking $V=\oline Y$ above we can therefore find
a cochain $\al=\sum_jc_j\al_j\in CF^0(\oline Y)$,
where the $\al_j$'s are generators and the coefficients are rational, such that
\be{equation}\label{cycle-cond}
\del'+d\al=0.
\end{equation}

For any non-positive integer $k$ set
\[\cm_k=\bigcup_{t\in\rbp}(M_t(E,k)\times\{t\}).\]
Let $v\in H^2(W;\z)/\torsion$ be any of the classes in
Lemma~\ref{Tred}~(ii). The zeros of the map
\[f_v:\R^{b^+}\to\R^{b^+}\]
defined in \Ref{fvdef} are precisely the parameter values of $\bt$
for which the $L^2$ $g(\bt)$-harmonic form representing $v$ is 
anti-self-dual. By making a small, $t$-independent perturbation
to the family of metrics $g(\bt)$ in some ball in $W$ where $g(t)$ is
independent of $t$ we can arrange that $0$
is a regular value of $f_v$ for each $v$.
Then $\cm_0$ will contain only finitely many reducible points.
After making a local perturbation to the ASD equations near each 
reducible point we can then arrange that all reducibles points in $\cm_0$
are regular (see Subsection~\ref{loc-str}).

Let $\cm_0'\subset\cm_0$ be the result of removing from $\cm_0$
a small neighbourhood about each reducible point, such that
$\prtl\cm_0'$ is a disjoint union of complex projective spaces.
For $k<0$ set $\cm_k'=\cm_k$.
(Recall that there are no reducibles in $\cm_k$ when $k<0$.)
Let $m$ be as in Theorem~\ref{h_ex} and set $n=(|w^2|-m)/2-1$. Since
\[\dim\,M_t(E,k)=8k-2c_1(E)^2-3(b^++1)\]
and $c_1(E)^2=w^2-b^+$, we have
\[\dim\,\cm_k=8k+4n+2m+1.\]

Consider the formal linear combination of moduli spaces
\[\tcm_k=\cm_k'+\sum_jc_j\cm_{\rho_j},\]
where $\rho_j$ is a projectively flat connection in $E|_Y$
representing $\al_j$ such that $\dim\,\cm_{\rho_j}=\dim\,\cm_k$.

Suppose the conclusion of Theorem~\ref{h_ex} does not hold, ie that
\be{equation}\label{ceil}
\lceil g/2\rceil+b^+-1<(|w^2|-m)/4.
\end{equation}
Set $q=\lceil g/2\rceil+b^+-1$. Then \Ref{ceil} is equivalent to
$2q\le n$.
Let $a\in H_2(X;\z)$ be as in the theorem and write
\[(x_1^2-4)^{\lceil g/2\rceil}\textstyle{\prod}_{i\ge2}(x_i^2-4))
=\sum_dP_d,\]
where $P_d$ is a homogeneous polynomial in $x_1,\dots,x_{b^+}$ of degree $2d$.
Cutting down as in the previous section we define
\[\hcm=\sum_{d=0}^q(a^mx^{n-2q}P_d)\cdot\tcm_{d-q},\]
which we regard as a formal linear combination of
smooth, oriented $1$-manifolds with boundary.

In the last two sections we will show that the number of boundary points of
$\hcm$,  counted with multiplicity, is non-zero, while the 
number of ends is zero. This contradiction will prove 
Theorem~\ref{h_ex}.

\section{Boundary points of $\hcm$}

Let $\ct\subset H^2(W;\z)$ be the torsion subgroup (which we may identify
with the torsion subgroup of $H^2(X;\z)$).

\be{prop}\label{red-count}
\[\#\prtl\hcm=\pm2^{b^+-1-2n-m}\,|\ct|\,\eta(\cl,w,a,m),\]
which is non-zero by the hypotheses of Theorem~\ref{h_ex}.
\end{prop}
\proof The boundary points of $\hcm$ all lie in $\cm'_0$, so
\[\#\prtl\hcm=\#(a^mx^{n-2q}P_q\cdot\prtl\cm'_0)
=\#(a^mx^n\cdot\prtl\cm'_0).\]

Now let $Q$ be the component of $\prtl\cm_0'$ corresponding to
some reducible point $D=([A],t)$. 
Fix the ordering of the corresponding splitting $E=L_1\oplus L_2$
and set $v=c_1(L_1)-c_1(L_2)$.
If we cut down $Q$ according to the monomial
$a^mx^n$ and count points with signs using the complex orientation of
$Q$ then the result is
\[\#(a^mx^n\cdot Q)=\eps2^{-2n-m}(v\cdot a)^m,\]
where $\eps=\pm1$ depends only on $m,n$, see \cite[Section~4]{Fr3}.
If $c=c_1(E)$ then the complex orientation of
$Q$ compare with the boundary orientation inherited from $\cm_0'$ by
\be{equation*}
\eps'\eps_{v,t}(-1)^{((v+c)/2)^2},
\end{equation*}
where $\eps'=\pm1$ is independent of $D$ and the ordering of $L_1,L_2$, and
$\eps_{v,t}=\pm1$ is the sign of the determinant of $df_v(t)$, see
Proposition~\ref{algebra}.
The last factor is taken from \cite[Prop.~3.25]{D2} and accounts
for whether the complex orientation of $\det(\del_D)$ agrees with the
orientation of the whole determinant line $\det(\del)$.

As explained in 
\cite[Section~4]{Fr3}, there are precisely $|\cal T|$
reducible points in $\cm_0$ for every pair $(\{v,-v\},t)$ where
$v\in H^2(W;\R)$ is one of the classes in Lemma~\ref{Tred} and $f_v(t)=0$.
Summing up we find that
\be{align*}
\#\prtl(a^mx^n\cdot\cm_0')
&=\pm{\textstyle\frac12}|\ct|
\sum_v\sum_{t\in f_v\inv(0)}\eps_{v,t}(-1)^{((v+c)/2)^2}
2^{-2n-m}(v\cdot a)^m\\
&=\pm|\ct|2^{-2n-m-1}\sum_v\text{deg}(v)
(-1)^{((v+c)/2)^2}(v\cdot a)^m.
\end{align*}
If $v=z-\sum_i\eps_is_i$ as in Lemma~\ref{Tred} then
$\deg(v)=\eps_1\cdots\eps_{b^+}$, and a simple
computation shows that
\[(-1)^{((z+w)/2)^2}=(-1)^{((v+c)/2)^2}\eps_1\dots\eps_{b^+}.\]
From this the Proposition follows immediately.\square

The reader may wish to check that 
$\eps_{v,t}(-1)^{((v+c)/2)^2}(v\cdot a)^m$ does indeed remain
unchanged when $v$ is replaced by $-v$.

\section{Ends of $\hcm$}\label{ns-ends}

There are two kinds of ends in $\hcm$. The first
kind arises from factorizations through flat connections over
$\rpy$. Because of our choice of ``limiting cycle''
over the end $\rpy$ (ie
condition~\Ref{cycle-cond}), the number of ends of this kind (counted with 
multiplicity) is zero. Indeed, by gluing theory (see \cite{D5}) there is
a finite collection of projectively flat connections $z_i$ in $E|_Y$ and for
each $i$ a finite subset $K_i\subset\cm_{z_i}$ equipped with a function
$K_i\to\z$ indicating multiplicities, such that the
ends of $\hcm$ can be identified with
\[\sum_iK_i\times\left(\check M(z_i,\theta) +
\sum_jc_j\check M(z_i,\rho_j)\right).\]
Therefore,
\[\#\{\text{ends of $\hcm$}\}=\sum_i(\#K_i)\,\la[z_i],\del'+\prtl\al\ra=0,\]
where $\al,\del'\in CF_4(Y)=CF^1(\oline Y)$.

The second type of ends in $\hcm$ arises from neck-stretching. More precisely,
these ends correspond to sequences
$\{(\bt(\nu),[A(\nu)])\}_{\nu=1,2,\dots}$ in $\hcm$ where $|\bt(\nu)|\to\infty$
as $\nu\to\infty$.
The goal of the remainder of this section is to prove that the number of ends
of this kind, counted with multiplicities, is zero.
We will make use of a recent result by Mu\~noz \cite{Mun1} which we
first explain.

Let $\Si$ be a closed surface of genus $g$ and
$F\to\Si$ the non-trivial
$\SO3$ bundle. Consider the affinely $\z/8$ graded
Floer cohomology group $\hf^*_g$ of the
$\SO3$ bundle $\bff=S^1\times F\to S^1\times\Si$.
Choosing an extension $\bff'\to D^2\times\Si$ of $\bff$ 
we can fix a $\z/8$-grading of $\hf^*_g$ by decreeing
that the element $1\in \hf^*_g$
defined by counting points in zero-dimensional moduli spaces in $\bff'$
has degree~$0$. (This grading is compatible with the canonical $\z/2$ grading
defined in \cite[Subsection~2.2]{Fr3}.)

The pair-of-pants cobordism gives rise to a product
\[\hf_g^p\otimes \hf_g^q\to \hf_g^{p+q}\]
which makes $\hf^*_g$ a commutative ring with unit.
Let $\tau$ be the degree~$4$ involution of $\hf^*_g$ defined by the class
$[S^1]\in H^1(S^1\times\Si;\z/2)$. This involution respects the product
in the sense that $x\cdot\tau(y)=\tau(x\cdot y)$ for
$x,y\in \hf^*_g$. Therefore the ring structure descends to the $\z/4$-graded
quotient $\thf^*_g=\hf^*_g/\tau$.

Let
\[\Psi_g:\text{Sym}(H_{\text{even}}(\Si))\otimes\La(H_1(\Si))\to \hf^*_g\]
be the invariant defined by the bundle $\bff'$, as explained in \cite{Fr3}.
Let $x\in H_0(\Si;\z)$ be the point class and
$\ga_1,\dots,\ga_{2g}\in H_1(\Si;\z)$
a symplectic basis with $\ga_j\cdot\ga_{j+g}=1$. Set $\al=2\Psi_g([\Si])$,
$\beta=-4\Psi_g(x)$, and $\ga=-2\sum_{j=1}^g\Psi_g(\ga_j\ga_{j+g})$, and let
$\ti\al,\ti\beta,\ti\ga$ denote the images of these classes in $\thf^*_g$.

The mapping class group of $\Si$ acts on $\hf^*_g$ and $\thf^*_g$ 
in a natural way. Mu\~noz shows that
the invariant part of $\thf^*_g$ is generated as a ring
by $\ti\al,\ti\beta,\ti\ga$, and he gives a
recursive description of the ideal of relations. The result we will use here
is the following.

\be{prop}[Mu\~noz]\label{munprop}\
\be{description}
\item[(i)]$\prod_{j=1}^g(\ti\beta+(-1)^j8)\in\ti\ga\thf^2_g$
\item[(ii)]$(\beta^2-64)^{\lceil g/2\rceil}\in\ga \hf^2_g$.
\end{description}
\end{prop}

Part~(ii) follows from part~(i), which
is stated explicitly in the proof of \cite[Proposition~20]{Mun1}.

We will now show that the number of ends of $\hcm$ coming from 
neck-stretching is zero.

Let $(\bt_\nu,[A_\nu])$, $\nu=1,2,\dots$ be a sequence of points in $\cm_r$
with $|\bt_\nu|\to\infty$ as $\nu\to\infty$. The transversality assumptions
imply that, after passing to a subsequence, $t_{i,\nu}$ will stay bounded
for all but one value of $i$.

Now fix $i$ and $\pm t_i=\tau\gg0$, so that $W$ contains a large cylinder
isometric to $[0,\pm t_i]\times S^1\times\Si^\pm_i$,
while the other $t_j$'s vary 
freely. Then the corresponding part $\hcm_{\pm t_i=\tau}$ of $\hcm$ is a finite
number of points (with multiplicities) which by gluing theory (see \cite{D5})
can be identified with a product
of instantons over $\R^2\times\Si^\pm_i$ and over $W\setminus\Si^\pm_i$
(with tubular end metrics). In particular, if we set
$E^\pm_i=E|_{S^1\times\Si^\pm_i}$, which is the pull-back of the
$\U2$ bundle over $\Si^\pm_i$ with $c_1=-1$, then
\[\#\hcm_{\pm t_i=\tau}=\phi^\pm_i\cdot\psi^\pm_i\]
for certain $\phi^\pm_i\in\cf^*(E^\pm_i)$ (measuring instantons over 
$\R^2\times\Si^\pm_i$) and $\psi^\pm_i\in\cf_*(E^\pm_i)$. 

If $i\ge2$ then
$\Si^\pm_i$ has genus~$1$. Thus $\cf^*(E^\pm_i)$ has only two generators,
in degrees differing by $4$ (see Appendix~\ref{u2conn}), and
\[\phi^\pm_i=\beta^2-64=0,\qquad\text{$i\ge2$}\]
{\em on chain level}. This is essential, because it implies not only that
$\#\hcm_{\pm t_i=\tau}=0$, but also that
$\psi^\pm_1$ is closed. The point here is that to prove $d\psi^\pm_1=0$
one must consider (suitably cut down) moduli spaces over
$W\setminus\Si^\pm_1$
of dimension~$1$, and these moduli spaces may have ends where
$t_j\to\pm\infty$ for some $j\ge2$.

Now let $\Phi^\pm,\Psi^\pm$ be the Floer (co)homology
classes of $\phi^\pm_1,\psi^\pm_1$, respectively.
By Mu\~noz' result we can write
\[\Phi^\pm=(\beta^2-64)^{\lceil g/2\rceil}=\ga e\]
for some $e\in\hf^*_g$, so
\[\#\hcm_{\pm t_1=\tau}=\ga e\cdot\Psi^\pm=e\cdot\ga\Psi^\pm.\]
But $\ga\Psi^\pm$ is a linear combination of Floer homology classes each of 
which can be defined by counting points in moduli spaces over
$W\setminus\Si^\pm_1$ cut down according to a monomial of the kind
$(\prod_{j=1}^Jz_j)(\prod_{i=2}^{b^+}x_i^{n_i})$
where $z_j\in H_{d_j}(W^-;\z)$ and $d_1=1$.
Since $H_1(W^-;\q)=0$ and we are using rational coefficients for the Floer
homology groups, this implies that $\ga\Psi^\pm=0$, as we will explain in the
next section. 

Given this, we conclude that if $\tau\gg0$ then
$\#\hcm_{\pm t_i=\tau}=0$ for all $i$, and Theorem~\ref{h_ex} follows.

\section{$\mu$-classes of loops}\label{mu-loop}

In this section we will give a simple proof of the vanishing result used
at the end of 
the previous section (ie $\ga\Psi^\pm=0$). This is essentially part of the
general assertion that the Donaldson-Floer invariants of $4$-manifolds
with boundary discussed in \cite[Subsection~2.3]{Fr3} are well-defined, in 
particular independent of the choice of geometric representatives for 
$\mu$-classes (see \cite{KM3}). One approach to this is to adapt the
proof for closed $4$-manifolds given in \cite{D3} and \cite{KM3},
but we prefer a more direct argument which also yields more precise 
information.

Because of the complexity of the proof of Theorem~\ref{h_ex} we will
consider the following simplified situation.
Let $X$ be an oriented Riemannian $4$-manifold with one tubular end
$\rpy$, and let $E\to X$ be a principal $\SO3$ bundle such that $E|_Y$ is
admissible in the sense of \cite{BD1}. To make sure moduli spaces are
orientable we assume 
$E$ lifts to a $\U2$ bundle, ie that $w_2(E)$ has an integral lift.

If $Y$ is a homology sphere we will assume there are no reducibles
in the instanton moduli spaces in $E$ considered below.

Set $X_0=X\setminus(\rpy)$ and let $\cb$ be the orbit space of $L^p_1$
connections in $E_0=E|_{X_0}$ modulo {\em even} gauge transformations 
of class $L^p_2$. Here $p$ should be an even integer greater than $4$,
and `even' means that the gauge transformation should lift to
$E\times_{\text{Ad}}\SU2$.
If $x_0\in X_0$ is a base-point then over the irreducible part
$\cb^*\subset\cb$ we have a natural
principal $\SO3$ bundle $\bp\to\cb^*$, the base-point fibration at $x_0$.
To this principal bundle and the adjoint representations of $\SO3$ on
itself, its Lie algebra, and its double cover $\Spin3$ we associate
three fibre bundles over $\cb^*$, which we denote
$\bg$, $\e$, and $\ti\bg$, respectively. Note that
$\bg$ is the bundle of fibrewise automorphisms of $\bp$.

Now consider a smooth loop $\lla:S^1\to X_0$ based at $x_0$. Holonomy 
along $\lla$ defines a smooth section $h$ of $\bg$.
For any flat connection $\al$ in $E|_Y$ set
\be{equation}\label{hol-fix}
Z_\al=\{[A]\in M(E,\al)\st h(A)=1\},
\end{equation}
where $M(E,\al)$ is the instanton moduli space in $E$ with flat limit $\al$.
After perturbing $h$ by a homotopy we may assume each $Z_\al$ is
transversely cut out, and we obtain a Floer cocycle
\[\phi=\phi(\lla)=\sum_\al(\#Z_\al)\al,\]
where the sum is taken over all equivalence classes of irreducible flat
connections $\al$ for which $\dim\,M(E,\al)=3$.
For rational coefficients the result we wish to prove is that the
cohomology class $[\phi]$ is a linear function of the homology class
$[\lla]$. However, we will actually prove a more precise statement, which
involves $\Spin3$ holonomy:

Up to isomorphism there are two spin structures on the pull-back bundle
$\lla^*E\to S^1$. Such a spin structure consists of a principal
$\Spin3$ bundle $\pi:Q\to S^1$ together with a bundle homomorphism
$Q\to\lla^*E$, see \cite{LM}. If $A$ is any connection in $E$ then there
is a unique connection $B$ in 
$Q$ such that $\pi(B)=\lla^*(A)$, and we can look at the holonomy of $B$. 
Because we are restricting to {\em even}
gauge transformations, this gives a section 
$\ti h$ of $\ti\bg$ which maps to $h$ under the covering $\ti\bg\to\bg$
and which depends on the choice of $Q$ only up to an overall sign.
Now, $\ti\bg$ is a bundle of Lie groups isomorphic to $\Spin3$, which we
can think of as the group of
unit quaternions, and the imaginary part $h_0$ of $\ti h$ is a 
section of the vector bundle $\e$.
Replacing the condition $h(A)=1$ in \Ref{hol-fix}
by $\ti h(A)=\pm1$ and $h_0(A)=0$, respectively, we get Floer cocycles
$\psi^\pm,\psi_0$, and these satisfy the relations
\[\phi=\psi_++\psi_-,\quad\psi_0=\psi_+-\psi_-.\]

Let $\Phi,\Psi^\pm,\Psi_0\in\hf^*(Y;\z)$ be the cohomology classes of
$\phi,\psi^\pm,\psi_0$, respectively. 
Note that $\Psi_0$ is independent of $\lla$, and since $\e$ is has odd
rank we have $2\Psi_0=0$.
Therefore,
\[\Psi_+=\Psi_-,\quad\Phi=2\Psi^\pm\qquad\text{modulo $2$-torsion}.\]
\be{prop}\
\be{description}
\item[(i)]The subset $\{\Psi_+(\lla),\Psi_-(\lla)\}\subset\hf^*(Y;\z)$
depends only on the class of $\lla$ in $H_1(X;\z)$.
\item[(ii)]$\lla\mapsto\Psi^\pm(\lla)$ defines a homomorphism
\[H_1(X;\z)\to\hf^*(Y;\z)/\text{$2$-torsion}.\]
\end{description}
\end{prop}
The Proposition is easily deduced from the following two lemmas:

\be{lemma}\label{psi-add}
If $\lla_1,\lla_2:S^1\to X_0$ are loops based at $x_0$ then for any 
compatible spin structures on $\lla_1^*E$, $\lla_2^*E$, and
$(\lla_1\circ\lla_2)^*E$, one has
\[\Psi_-(\lla_1\circ\lla_2)=\Psi_-(\lla_1)+\Psi_-(\lla_2).\]
\end{lemma}
{\em Proof of Lemma:} Set $G=\Spin3$ and let $S\bh_0$ denote the unit sphere
in the space $\bh_0$ of pure quaternions. (We identify $\bh_0$ with the
Lie algebra of $G$.)
Define a subset
$V\subset G\times G$ by
\[V=\{(\exp(sx),\exp(tx))\st x\in S\bh_0,\; s,t\in[0,\pi],\;s+t\ge\pi\}.\]
If we ignore the singular points $(1,-1),(-1,1),(-1,-1)$, then $V$ is a 
smooth, orientable $4$-manifold with boundary, and
\[\prtl V=(\{-1\}\times G)\cup(G\times\{-1\})\cup D,\]
where $D=\{(g_1,g_2)\in G\times G\st g_1g_2=-1\}$.

If $M$ is any closed, oriented $3$-manifold and $f_1,f_2:M\to G$ 
smooth maps such that $(f_1,f_2):M\to G\times G$ misses the three
singular points of $V$ and is transverse to both $V$ and $\prtl V$ then
\[0=\#\prtl(f_1\times f_2)\inv(V)=
\eps_1\deg(f_1)+\eps_2\deg(f_2)-\eps_3\deg(f_1\cdot f_2)\]
for some constants $\eps_1,\eps_2,\eps_3=\pm1$.
Taking $M=G$, $f_1=\text{Id}$, $f_2\equiv1$ gives $\eps_3=\eps_1$,
and similarly $\eps_3=\eps_2$.

Let $\ti{\bg}\otimes\ti{\bg}$ denote the fibrewise
product of $\ti{\bg}$ with
itself. Since $V$ is $\text{Ad}_G$-invariant it defines a subset
$W\subset\ti{\bg}\otimes\ti{\bg}$ which is a fibre bundle over $\cb^*$
with fibre $V$. Now let $\ti h_j$ be the section of $\ti\bg$ obtained from
$\lla_j$, and set 
\[Z'_\al=\{[A]\in M(E,\al)\st(\ti h_1(A),\ti h_2(A))\in W\}.\]
Then the Floer cochain
\[b=\uset{\dim M(E,\al)=2}\sum(\#Z'_\al)\al\]
satisfies
\[db=\eps_1\psi_-(\lla_1)+\eps_2\psi_-(\lla_2)-
\eps_3\psi_-(\lla_1\circ\lla_2)\]
with the same constants $\eps_j$ as above.\square

\be{lemma}If $\lla:S^1\to X_0$ represents the zero class in $H_1(X;\z)$
then there is a spin structure on $\lla^*E$ for which $\Psi_-(\lla)=0$.
\end{lemma}
\proof There is a smooth, compact,
oriented, connected $2$-manifold $\Si$ with one boundary component (which
we identify with $S^1$), together with a smooth map
$f:\Si\to X_0$ such that $f|_{\prtl\Si}=\lla$. Choose a spin structure on
$f^*E$. It is easy to see that if $\mu:S^1\to\Si$ is any 
loop based at $z\in\prtl\Si$ then $\Psi_-(f\circ\mu)$ depends only on
the class of $\mu$ in the fundamental group
$\pi_1(\Si,z)$. Since the identity map
$S^1\to\prtl\Si$ is a product of commutators in $\pi_1(\Si,z)$, it
follows from Lemma~\ref{psi-add} that $\Psi_-(\lla)=0$.
\square

\appendix

\section{Determinant line bundles}\label{sec:detlb}

This appendix gives an account of the construction
of determinant line bundles for families of Fredholm operators, taking
care of some signs that many authors have overlooked.

Determinant line bundles for families of elliptic operators arise naturally
in gauge theory and symplectic geometry
in connection with orientations of moduli spaces,
see \cite{DK,FH}.
In general, if $\{T(x)\}_{x\in C}$ is a continuous
family of Fredholm operators between two Banach spaces then
the determinant line bundle $\det(T)$, as a set, is the disjoint union of all
the determinants $\det(T(x))$ as $x$ varies through $C$, see below. 
There is a natural collection of local trivializations of $\det(T)$ which one
can attempt to use to make $\det(T)$ a topological line bundle over $C$.
What many authors seem to have overlooked is that the overlap
transformations are not in general continuous. We resolve this problem by
dividing the set of local trivializations into two parts, thereby obtaining
an ``even'' and an ``odd'' topology on $\det(T)$. The two topologies
are interchanged by the involution of $\det(T)$ which is multiplication
with $(-1)^{\dim(\ker(T(x)))}$ on the fibre $\det(T(x))$.

For the sake of simplicity, we will in the main text
only consider the {\em odd} topology (this prevents a sign in 
Equation~\ref{gadd}).

\subsection{Determinant lines}\label{subsec:detl}

We first study the determinant line of a single Fredholm operator.

For a finite dimensional vector space $A$ we denote by $\det(A)$ the
highest exterior power of $A$. If
\[0\to A_0\oset{\al_0}\to A_1\oset{\al_1}\to A_2\to0\]
is an exact sequence of linear maps between finite dimensional vector
spaces then there is a natural isomorphism
\be{equation}\label{sign-convention}
\det(A_0)\otimes\det(A_2)\oset\approx\to\det(A_1),\quad
x_0\otimes x_2\mapsto\al_0(x_0)\wedge s(x_2),
\end{equation}
where $s:A_2\to A_1$ is any right inverse of $\al_1$. More generally, if
\[0\to A_0\oset{\al_0}\to A_1\oset{\al_1}\to\dots\oset{\al_{r-1}}\to A_r\to0\]
is an exact sequence of linear maps one gets a natural isomorphism
\[\uset{\text{$i$ even}}\otimes\det(A_i)\oset\approx\to
\uset{\text{$i$ odd}}\otimes\det(A_i).\]
Now let $V,W$ be Banach spaces over $\f=\R$ or $\co$,
and $S:V\to W$ a Fredholm operator. 
The determinant line of $S$ is by definition
\[\det(S)=\det(\ker(S))\otimes(\det(\coker(S)))^*.\]
If $f:\f^n\to W$ is a linear map set
\[S_f:V\oplus\f^n\to W\oplus\f^n,\quad(v,z)\mapsto(Sv+fz,0),\]
and let $S\oplus f:V\oplus\f^n\to W$ be the $W$-component of $S_f$.
Then there is a natural exact sequence
\[0\to\ker(S)\to\ker(S_f)\to\f^n\oset{\bar f}\to\coker(S)\to\coker(S_f)
\to\f^n\to0,\]
where $\bar f$ is the map $f$ followed by the projection onto $\coker(S)$.
This sequence gives rise to a natural isomorphism
\be{equation}\label{det-iso}
\det(S)\oset\approx\to\det(S_f).
\end{equation}

Next we consider a pair of linear maps $f_j:\f^{n_j}\to W$, $j=1,2$. Set
\[f_1\oplus f_2:\f^{n_1+n_2}=\f^{n_1}\oplus\f^{n_2}\to W,\quad
(z_1,z_2)\mapsto f_1z_1+f_2z_2.\]
Then we have a diagramme of isomorphisms
\be{equation}\label{ucdiagr}
\be{array}{ccc}
\dsp\det(S) & \oset{\textstyle\lla_1}\longrightarrow
 & \det(S_{f_1})\\
\dsp\lla_3\downarrow\hphantom{\lla_3} & & \hphantom{\lla_2}\downarrow\lla_2\\
\det(S_{f_1\oplus f_2}) & =  & \det((S_{f_1})_{f_2})
\end{array}
\end{equation}
where all maps except the bottom horizontal one are instances of \Ref{det-iso}.
Consider the three maps
\be{equation}\label{threemaps}
\f^{n_1}\to\coker(S),\quad\f^{n_2}\to\coker(S\oplus f_1),
\quad\f^{n_1+n_2}\to\coker(S)
\end{equation}
obtained from $f_1,f_2,f_1\oplus f_2$ in the obvious way. 
Let $K_1,K_2,K_{12}$ be the kernels of the maps in \Ref{threemaps}, and
let $L_1,L_2,L_{12}$ be
the images of the same maps.
\be{lemma}\label{square-sign}
$\lla_3=(-1)^{\dim(L_1)\dim(K_2)}\lla_2\lla_1$.
\end{lemma}
We remark that the sign in the Lemma does not depend on the sign-convention 
in \Ref{sign-convention} (ie the order of $\al_0(x_0),s(x_2)$).

\proof The proof is an exercise in understanding the definitions involved.
We will merely indicate where the sign comes from. Note that
there is a commutative diagramme
\be{equation}\label{kkll}
\be{array}{ccccccccc}
\dsp
& & 0 & & 0 & & 0 & & \\
\dsp
& & \downarrow & & \downarrow & & \downarrow & & \\
\dsp
0 & \to & K_1 & \to & \f^{n_1} & \to &
L_1 & \to & 0 \\
\dsp
& & \downarrow & & \downarrow & & \downarrow & & \\
\dsp
0 & \to & K_{12} & \to & \f^{n_1+n_2} & \to &
L_{12} & \to & 0 \\
\dsp
& & \downarrow & & \downarrow & & \downarrow & & \\
\dsp
0 & \to & K_2 & \to & \f^{n_2} & \to &
L_2 & \to & 0 \\
\dsp
& & \downarrow & & \downarrow & & \downarrow & & \\
& & 0 & & 0 & & 0 & & 
\end{array}
\end{equation}
where all rows and columns are exact. This defines two isomorphisms
\[\det(K_1)\otimes\det(L_1)\otimes\det(K_2)\otimes\det(L_2)
\oset\approx\to\det(\f^{n_1+n_2}),\]
which differ by the factor $(-1)^{\dim(L_1)\dim(K_2)}$.
One isomorphism uses the top and bottom horizontal exact sequences and
then the middle vertical sequence. The other isomorphism uses
the remaining three exact sequences in the diagramme.

The point is that when computing $\lla_2\lla_1$
one uses the first isomorphism, while
$\lla_3$ involves the other one.
The remaining details are left to the reader.\square

\subsection{Determinant line bundles}

Let $\cb(V,W)$ the Banach space of bounded operators
from $V$ to $W$, and $\text{Fred}(V,W)\subset\cb(V,W)$ the open subset
consisting of all Fredholm operators. If $C$ is a space and 
$T:C\to\text{Fred}(V,W)$ a continuous map we define $\det(T)$ as a set by
\[\det(T)=\bigcup_{x\in C}\{x\}\times\det(T(x)).\]
Note that if $T(x)$ is surjective for every $x$ then 
$\ker(T)\subset C\times V$, the union of all sets $\{x\}\times\ker(T(x))$, 
is a topological vector bundle over $C$,
so in this case $\det(T)$ has a natural topology. For a general $T$, 
we will topologize $\det(T)$ by essentially
specifying a set of local trivializations.
Since the surjective operators are open in $\cb(V,W)$, every point in $C$ has
an open neighbourhood $U$ for which there exists a linear map
$f:\f^n\to W$ such that
$T(x)\oplus f:V\oplus\f^n\to W$ is surjective for every $x\in U$.
By definition,
\[\det(T_f(x))=\det(T(x)\oplus f)\otimes\det(\f^n)^*,\]
hence $\det(T_f)|_U$ has a natural topology.
Moreover, fibrewise application of \Ref{det-iso} gives a bijection
\[\lla_{U,f,n}:\det(T)|_U\to\det(T_f)|_U.\]
\be{prop}
\be{description}
\item[(i)]If $\eps\in\z/2$ then $\det(T)$ has a unique topology such
that the projection $\det(T)\to C$ is continuous and all 
maps $\lla_{U,f,n}$ with $n\equiv\eps\mod2$ are 
homeomorphisms. The corresponding space $\det_\eps(T)$ is a topological
line bundle over $C$.
\item[(ii)]The fibre preserving map
\[\dett_0(T)\to\dett_1(T)\]
which is multiplication by $(-1)^{\dim(\ker(T(x)))}$
on the the fibre $\det(T(x))$,
is a homeomorphism.
\end{description}
\end{prop}

In the main text we will, as already mentioned, take $\eps=1$.

\proof The only remaining issue is continuity of overlap transformations.
Suppose $f_j:\f^{n_j}\to W$, $j=1,2$ are linear maps such that
$T(x)\oplus f_j$ is surjective for all $x$. Then we have a diagramme of
bijective maps
\be{equation*}
\be{array}{ccccccc}
\dsp\det(T_{f_1}) & \oset{\textstyle\lla_1}\longleftarrow &
\det(T) & = & \det(T) & \oset{\textstyle\lla_2}\longrightarrow &
\det(T_{f_2})\\ \dsp
\beta_1\downarrow\hphantom{\beta_1} & & \downarrow & & \downarrow & & 
\hphantom{\beta_2}\downarrow\beta_2\\
\det((T_{f_1})_{f_2}) & = & \det(T_{f_1\oplus f_2}) & = & 
\det(T_{f_2\oplus f_1}) & = & \det((T_{f_2})_{f_1})
\end{array}
\end{equation*}
where all maps marked with an arrow are instances of \Ref{det-iso}.
The middle square obviously commutes. In the fibres over $x\in C$
the square to the left commutes up to multiplication with
$(-1)^{n_2\dim(\coker(T(x)))}$, by Lemma~\ref{square-sign}, and similarly
the square to the right commutes up to multiplication with
$(-1)^{n_1\dim(\coker(T(x)))}$. It follows that
\[\lla_2\lla_1\inv=(-1)^{(n_1+n_2)\dim(\coker(T(x)))}\beta_2\inv\beta_1.\]
Since $\beta_2\inv\beta_1$ is continuous, and
\[\text{index}(T(x))=\dim(\ker(T(x)))-\dim(\coker(T(x)))\]
is locally constant, the Proposition follows.\square

\subsection{Isomorphisms between determinant line bundles}
\label{cont-det}

Suppose $g:C\to\cb(\f^p,W)$ is continuous and $\eps_1,\eps_2\in\z/2$.
We will define a canonical isomorphism of line bundles
\[\ga:\dett_{\eps_1}(T)\to\dett_{\eps_2}(T_g).\]
If $n_j\equiv\eps_j\pod2$,
and $f_j:\f^{n_j}\to W$ are linear maps such that $T\oplus f_1$ and
$T\oplus g\oplus f_2$ are both surjective in some open set
$U\subset C$, then $\det(T_{f_1})$ and $\det(T_{g\oplus f_2})$
can both be identified with
$\det(T_{g\oplus f_1\oplus f_2})$ over $U$. This gives an isomorphism
$\det(T)|_U\to\det(T_g)|_U$, which is easily seen to be 
independent of $f_1,f_2$, hence patch together to give
the desired isomorphism $\ga$.

In explicit computations in may be useful to have $\ga$ expressed in terms
of the bijective map
\[\ga':\dett_{\eps_1}(T)\to\dett_{\eps_2}(T_g)\]
defined in \Ref{det-iso}. Set
\[d_1=\dim(\coker(T)),\quad d_2=\dim(\coker(T\oplus g)),\]
regarded as functions on $C$. Then Lemma~\ref{square-sign} and a simple
computation gives:
\be{prop}\label{gagap}
$\ga=(-1)^{d_1\eps_1+d_2\eps_2+d_1p+d_1d_2+d_2}\ga'$.
\end{prop}

Finally, we remark that the isomorphism $\ga$ is functorial in the sense
that if $g_j:C\to\cb(\f^{p_j},W)$ is a continuous map for $j=1,2$ and if
$\eps_1,\eps_2,\eps_3\in\z/2$ then the composition of the 
isomorphisms
\[\dett_{\eps_1}(T)\to\dett_{\eps_2}(T_{g_1})\to
\dett_{\eps_3}(T_{g_1\oplus g_2})\]
agrees with the isomorphism
$\det_{\eps_1}(T)\to\det_{\eps_3}(T_{g_1\oplus g_2})$.

\section{On projectively flat $\U2$ connections over $3$-manifolds}
\label{u2conn}

In this appendix we will prove a simple result which describes (modulo 
$\SU2$ gauge equivalence)
projectively flat $\U2$ connections over a closed, oriented
$3$-manifold in terms of flat $\SU2$ connections over the complement
of a suitable link. The existence of
such a correspondence was pointed out in \cite[p198]{BD1}. We will then
apply this result to compute the Floer chain complex of a non-trivial 
$\SO3$ bundle over the $3$-torus.

Let $Y$ be a closed, oriented $3$-manifold, $E\to Y$ a rank~$2$
Hermitian vector bundle, and $L=\La^2E$ the determinant line bundle of $E$.
Let $s$ be a regular, smooth section of $L$. The zero-set of $s$ is then
a disjoint union $\ga=\coprod_j\ga_j$ of embedded circles. 
Let $N\approx D^2\times\ga$ be a closed tubular neighbourhood of $\ga$, and
let $V\subset Y$ denote the complement of
$(\frac12\text{int}D^2)\times\ga$. By modifying $s$ if necessary we may
arrange that $|s|=1$ in $V$. Let $\hat A$ be a unitary connection in $L$ such
that $\nabla_{\hat A}s=0$ in $V$. 
Since $L|_N$ is trivial there is a unitary isomorphism
$L|_N\approx L_0\otimes L_0$ where $L_0=N\times\co$. Let $\hat A_0$ be the
connection in $L_0$ induced by $\hat A$. Of course, $F(\hat A)=2F(\hat A_0)$
in $N$. 

For each $j$ choose a point $z_j\in\ga_j$ and define $m_j=S^1\times\{z_j\}$.
Then
\[\hol_{m_j}(\hat A_0)=\exp\left(-\int_{D^2\times z_j}F(\hat A_0)\right)=-1,\]
where the second equality is a basic fact from Chern-Weil theory.

Let $M_1$ denote the moduli space of projectively flat unitary connections in
$E$ which induces $\hat A$ in $L$, modulo automorphisms of $E$ of
determinant~$1$. If $P$ is the $\SO3$ bundle associated to $E$ then $M_1$
can be identified with the moduli space $\fl_S(P)$
of flat connections in $P$ modulo
{\em even} automorphisms of $P$, ie those that lift to
$P\times_{\text{Ad}(\SO3)}\SU2$. 
If $[A]\in M_1$ then $A$ and
$\hat A_0\oplus\hat A_0$ restrict to gauge
equivalent connections over $D^2\times z_j$, hence $\hol_{m_j}(A)=-1$.

Let $M_2$ denote the moduli space of flat unitary connections in
$E|_V$ satisfying $\nabla_As=0$ and $\hol_{m_j}(A)=-1$ for all
$j$, modulo automorphisms of $E|_V$ of determinant~$1$.

\be{prop}\label{m12}
The restriction map $r:M_1\to M_2$ is a bijection.
\end{prop}

\proof
Given $[B]\in M_2$ it is easy to find a flat $\SU2$ connection $A_0$ over $N$
such that if $A=A_0\otimes\hat A_0$ then $A$ and $B$ restrict to gauge
equivalent $\U2$ connections over $N\cap V$. Gluing these together along
$N\cap V$ produces an element $\al\in M_1$ such that $r(\al)=[B]$.
Hence $r$ is surjective.
It is easy to see that $r$ is also injective.\square

\be{prop}If $P\to T^3$ is any non-trivial $\SO3$ bundle then $\fl_S(P)$ 
has exactly two elements. These are both non-degenerate and differ in index 
by $4$.
\end{prop}

\proof Recall that $\SO3$ bundles over a compact $3$-manifold (or over any
finite CW-complex of dimension $\le3$) are determined up to isomorphism
by their second Stiefel-Whitney class. Choose an indivisible class
$c\in H^2(T^3;\z)$ which is a lift of $w_2(P)$. Then there exists a 
diffeomorphism of $T^3$ (defined by some element of $\SL(3,\z)$) such that
$f^*c$ is the Poincar\'e dual of $[S^1]\times1\times1$. We may therefore
assume that $P=S^1\times P_0$, where $P_0\to T^2$ is the non-trivial 
$\SO3$ bundle.
Applying Proposition~\ref{m12} with $\ga=S^1\times\text{pt}$ we find that
$\fl_S(P)$ has exactly two elements $\al_\pm$, given by the
representations of
$\pi_1(S^1\times(T^2\setminus\text{pt}))$ into $\SU2=\Sp1$
which take the three standard generators 
to $\pm1,i,j$, respectively, where $1,i,j,k$ is the usual
basis for the quaternion algebra. 

It is not hard to see that the corresponding flat connection in $P$ is
non-degenerate, by observing that it pulls back to the trivial connection
under the $4$-fold covering $T^3\to T^3$, $(r,s,t)\mapsto(r,s^2,t^2)$.
We know apriori that there is 
a degree~$4$ involution of $\fl_S(P)$ (see \cite{Fr3}), so the index
difference of $\al_\pm$ must be $4\mod8$.
\square

\textsc{Institut des Hautes \'Etudes Scientifiques,\\
F-91440 Bures-sur-Yvette, France}

\end{document}